\documentclass[12pt]{article}
\usepackage{amssymb,amsmath,mathrsfs}
\usepackage{hyperref}
\usepackage{graphicx}
\usepackage{color}
\usepackage[OT2,OT1]{fontenc}
\usepackage{upquote}
\usepackage{authblk}
\usepackage[numbers,sort&compress]{natbib}

\begin{document}

\title{\LARGE\bf A method to reduce the Lehmer measure in a multi-term Machin-like formula for $\pi$}

\bigskip
\author[1, 2]{\small Sanjar M. Abrarov}
\author[2, 3, 4]{\small Rehan Siddiqui}
\author[3, 4]{\small Rajinder K. Jagpal}
\author[1, 2, 4]{\small \\ Brendan M. Quine}

\affil[1]{\scriptsize Thoth Technology Inc., Algonquin Radio Observatory, Achray Road, RR6, Pembroke,~ON,~Canada,~K8A~6W7 \normalsize}
\affil[2]{\scriptsize Dept. Earth and Space Science and Engineering, York University, 4700 Keele St., Toronto,~ON,~Canada,~M3J~1P3 \normalsize}
\affil[3]{\scriptsize Epic College of Technology, 5670 McAdam Rd., Mississauga, ON, Canada, L4Z 1T2 \normalsize}
\affil[4]{\scriptsize Dept. Physics and Astronomy, York University, 4700 Keele St., Toronto, ON, Canada, M3J 1P3 \normalsize}

\date{October 19, 2022}
\maketitle
\maketitle

\begin{abstract}
Previously we have proposed a new method of transforming quotients into integer reciprocals in the Machin-like formulas for $\pi$. As a further development, here we show how to generate a multi-term Machin-like formula for $\pi$ with a reduced Lehmer measure. The Mathematica codes validating these results are presented.
\vspace{0.25cm}
\vspace{0.25cm}
\end{abstract}

\section{Introduction}

The Machin-like formulas for $\pi$ can be represented in the following form \cite{Lehmer1938, Abeles1993,Beckmann1971,Berggren2004,Borwein2008}
\begin{equation}\label{eq_1}
\frac{\pi}{4} = \sum\limits_{j = 1}^J A_j\arctan\left(\frac{1}{B_j}\right), \qquad {A_j},{B_j} \in \mathbb{Q}.
\end{equation}
In many known equations the values ${B_j}$ are integers. Among many other formulas for $\pi$ \cite{Beckmann1971,Berggren2004,Borwein2008,Fitzhugh2013,Hwang2017}, the identities of kind \eqref{eq_1} involving arctangent function represent a particular interest due to their potential for rapid convergence.

Since the Maclaurin series expansion of the arctangent function is given by
\begin{equation}\label{eq_2}
\arctan\left(x\right) = x - \frac{x^3}{3} + \frac{x^5}{5} - \frac{x^7}{7} + \cdots = \sum\limits_{n = 0}^\infty \frac{\left(-1\right)^n x^{2n + 1}}{2n + 1},\quad - 1\leqslant x\leqslant 1,
\end{equation}
we can conclude that in order to improve convergence rate the numbers ${B_j}$ in equation \eqref{eq_1} should be as large as possible by absolute values. Furthermore, the computational efficiency of a given Machin-like formula \eqref{eq_1} for $\pi$ is obviously higher if the number of the terms $J$ is smaller. It should be noted that implementation of the binary splitting algorithm can also be used to accelerate computation of the series \eqref{eq_2} (see for example \cite{Brent2010}).

Lehmer introduced a measure \cite{Lehmer1938,Abeles1993,Tweddle1991,Wetherfield1996}
\begin{equation}\label{eq_3}
\mu  = \sum\limits_{j = 1}^J\frac{1}{\log_{10}\left(\left|B_j\right|\right)}.
\end{equation}
Nowadays, the equation \eqref{eq_3} is known as the Lehmer measure. We take absolute values $\left|B_j\right|$ since the constants $B_j$ in equation \eqref{eq_1} may generally be negative.

According to Lehmer, the measure \eqref{eq_3} shows "amount of labor" that is required for computation of $\pi$ for some given Machin-like formula. Specifically, the lower Lehmer measure $\mu$ indicates the higher computational efficiency of the Machin-like formula \eqref{eq_1} for $\pi$. For example, the original Machin formula for $\pi$ \cite{Lehmer1938,Abeles1993,Beckmann1971,Berggren2004,Borwein2008,Jansson2019,Nimbran2010,Chamberland2019}
\begin{equation}\label{eq_4}
\frac{\pi}{4} = 4\arctan\left(\frac{1}{5}\right) - \arctan\left(\frac{1}{239}\right)
\end{equation}
that was found by English mathematician John Machin in 1706 has the Lehmer measure $\mu  \approx 1.85113$. However, the following Machin-like formula for $\pi$
\begin{equation}\label{eq_5}
\frac{\pi}{4} = 12\arctan\left(\frac{1}{18}\right) + 8\arctan\left(\frac{1}{57}\right) - 5\arctan\left( \frac{1}{239}\right)
\end{equation}
that was discovered by Gauss has the Lehmer measure $\mu\approx 1.78661$ \cite{Lehmer1938,Abeles1993}. This signifies that the formula \eqref{eq_5} is more efficient in computation of $\pi$ than the formula \eqref{eq_4} as its Lehmer measure is smaller. More detailed descriptions and significance of the Lehmer measure \eqref{eq_3} in computation of $\pi$ can be found in literature \cite{Wetherfield1996}.

Previously we have shown that the following identity \cite{Abrarov2017a} (see also \cite{WolframCloud})
$$
\frac{\pi}{4} = 2^{k-1}\arctan\left(\frac{\sqrt{2 - a_k}}{a_k}\right),\qquad k\geqslant 1,
$$
where $a_k = \sqrt{2 + a_{k - 1}}$, $a_0 = 0$, can be rearranged as
\begin{equation}\label{eq_6}
\frac{\pi}{4} = 2^{k-1}\arctan\left(\frac{1}{\alpha}\right) + \arctan\left(\frac{1}{\beta}\right),
\end{equation}
where the constant $\alpha $ can be taken as
\begin{equation}\label{eq_7}
\alpha  = \left\lfloor\frac{a_k}{\sqrt{2 + a_k}}\right\rfloor
\end{equation}
for any integer $k$ greater than $1$.

Once $\alpha$ is known, the constant $\beta$ in equation \eqref{eq_6} can be found by using the following equation \cite{Abrarov2017a,WolframCloud}
\begin{equation}\label{eq_8}
\beta = \frac{2}{\left(\left(\alpha + 1\right)/\left(\alpha - 1\right)\right)^{2^{k - 1}} - i} - i.
\end{equation}
Application of equation \eqref{eq_8} is not optimal since increasing $k$ drastically slows down computation of the constant $\beta$. Fortunately, as we have shown in our work \cite{Abrarov2017b} the constant $\beta$ in equation \eqref{eq_6} can be found more efficiently by using a two-step iteration
\begin{equation}\label{eq_9}
\left\{ \begin{gathered}
\sigma_n = \sigma_{n - 1}^2 - \tau_{n - 1}^2\hfill \\
\tau_n = 2\sigma_{n - 1}^2\tau_{n - 1}^2,\qquad\qquad n = \left\{2,3,4,\ldots,k\right\},\hfill \\ 
\end{gathered}\right.
\end{equation}
with initial values defined as
\[
\label{eq_10a}\tag{10a}
\sigma_1 = \frac{\alpha^2 - 1}{\alpha^2 + 1}
\]
and
\[
\label{eq_10b}\tag{10b}
\tau_1 = \frac{2\alpha^2}{\alpha^2 + 1}
\]
such that
\[
\label{eq_10c}\tag{10c}
\beta = \frac{\sigma_k}{1 - \tau_k}.
\]

A method, known as the Todd's process, is commonly used to generate the Machin-like formulas for $\pi$  \cite{Wetherfield1996} (see also \cite{Todd1949}). However, this method is not simple and requires tedious and complicated manipulations with large-size matrices based on a set of some prime numbers \cite{Wetherfield1996}. Furthermore, it seems that computer-aided algorithms based on Todd's process still cannot generate the Machin-like formulas for $\pi$ containing only integer reciprocals with $\mu$ below $1$. For example, the following Machin-like formula for $\pi$ containing only integer reciprocals (see \cite{Wetherfield2004} and web-links therein)
\setcounter{equation}{10}
\begin{equation}\label{eq_11}
\begin{aligned}
\frac{\pi}{4}=&83\arctan\left(\frac{1}{107}\right) + 17\arctan\left(\frac{1}{1710}\right) - 22\arctan\left(\frac{1}{103697}\right)\\
&-24\arctan\left(\frac{1}{2513489}\right) - 44\arctan\left(\frac{1}{18280007883}\right)\\
&+12\arctan\left(\frac{1}{7939642926390344818}\right)\\
&+22\arctan\left(\frac{1}{3054211727257704725384731479018}\right)
\end{aligned}
\end{equation}
has the Lehmer measure $\mu\approx 1.34085$. To the best of our knowledge this is the smallest Lehmer measure ever reported for the Machin-like formulas for $\pi$ containing only integer reciprocals.

In our previous publication we have proposed a new method that transforms quotients into integer reciprocals in the Machin-like formulas for $\pi$ \cite{Abrarov2021a}. In particular, we have shown how this approach can be used to generate a Machin-like formula for $\pi$ that can retain without any quotients the Lehmer measure $\mu$ below one. As a further development of our work \cite{Abrarov2021a}, here we show how to generate the Machin-like formula for $\pi$ with the Lehmer measure $\mu$ less than $1$ even at a smaller values of the integer $k$. As an example we show that for $k=4$ the corresponding Machin-like formula \eqref{eq_22} for $\pi$ has the Lehmer measure $0.956916$. Therefore, this approach may be promising for efficient computation of digits of $\pi$.

\section{Preliminaries}

As an example, we can take  $k = 4$. Substituting this value into equation \eqref{eq_7} we have
$$
\alpha  = \left\lfloor\frac{\sqrt{2 + \sqrt{2 + \sqrt{2 + \sqrt 2}}}}{\sqrt{2 - \sqrt{2 + \sqrt{2 + \sqrt 2}}}}\right\rfloor = 10.
$$
Then, using equations \eqref{eq_10a}, \eqref{eq_10b} and \eqref{eq_10c} together with two-step iteration \eqref{eq_9}, we can find that
$$
\beta = -\frac{147153121}{1758719}.
$$
Consequently, substituting these values $\alpha $ and $\beta $ into equation \eqref{eq_6} yields
\begin{equation}\label{eq_12}
\frac{\pi}{4} = 8\arctan\left(\frac{1}{10}\right) - \arctan\left(\frac{1758719}{147153121}\right).
\end{equation}

We use can Mathematica to verify this two-term Machin-like formula for $\pi$ in two ways. The first way is to use a built-in algorithm of the Mathematica by equating logically (with == symbols) the left and right sides of equation \eqref{eq_12}. The following Mathematica code: \\[0.25cm]
\small
{\ttfamily{Pi/4==8*ArcTan[1/10]-ArcTan[1758719/147153121]}} \\[0.25cm]
\normalsize
\noindent validates the identity \eqref{eq_12} by returning {\ttfamily True}. Although this code confirms correctness of this identity, it is not clear how the built-in algorithm of Mathematica validates it.

Another way of validation that follows from equation \eqref{eq_1}, implies that the real and imaginary parts of the product (detailed description can be found in \cite{Jansson2019})
$$
\prod\limits_{j = 1}^J \left(B_j + i\right)^{A_j}\propto\left(1 + i\right)
$$
must be equal to each other. For example, there is an elementary proof for the original Machin formula \eqref{eq_4} for $\pi$ \cite{Guillera2009}
$$
(5 + i)^4(259 + i)^{-1}=2(1 + i).
$$
Therefore, the following Mathematica code:\\[0.25cm]
\small
{\ttfamily coeff=(10+I)\string^8*(147153121/1758719+I)\string^-1;\\
Re[coeff]==Im[coeff]}\\[0.25cm]
\normalsize
also validates the equation \eqref{eq_12} by returning {\ttfamily True} since the real and imaginary parts of the product
$$
\left(10 + i\right)^8\left(\frac{147153121}{1758719 + i}\right)^{-1} = \frac{1758719}{2}\left(1 + i\right)
$$
are equal to each other. Further, we will use the second way of validation since it is explicit and, therefore, intuitive for understanding.

The Machin-like formula \eqref{eq_12} for $\pi$ contains a quotient in the argument of the second arctangent function. This quotient is undesirable as its large numerator $1758719$ increases the argument of the arctangent function that slows down the convergence rate. Moreover, it also contributes for more digits that also slows down the computation when expansion series like \eqref{eq_2} for the arctangent function is used. Therefore, it is very desirable to look for the Machin-like formulas for $\pi$, where all constants ${B_j}$ are integers. In order to resolve this problem we proposed a new method based on the following identity (see \cite{Abrarov2021a} for derivation)
\begin{equation}\label{eq_13}
\arctan\left(\frac{1}{z}\right) = \arctan\left(\frac{1}{\left\lfloor z \right\rfloor}\right) + \arctan \left(\frac{\left\lfloor z \right\rfloor - z}{1 + z\left\lfloor z \right\rfloor}\right),\quad z\notin [0,1).
\end{equation}

Application of the identity \eqref{eq_13} has appeared to be very efficient. Consider as an example a remarkable Machin-like formula for $\pi$ that was found by Wetherfield in 2004 (see \cite{Wetherfield2004} and the web-links therein)
\begin{equation}\label{eq_14}
\begin{aligned}
\frac{\pi}{4} = &83 \arctan\left(\frac{1}{107}\right) + 17\arctan\left(\frac{1}{1710}\right) - 22\arctan\left(\frac{1}{103697}\right) \\
&-12\arctan\left(\frac{2}{2513489}\right) - 22\arctan\left(\frac{2}{18280007883}\right).
\end{aligned}
\end{equation}
The Lehmer measure of this formula is small $\mu \approx 1.26579$ (see definition \eqref{eq_3} for computation of $\mu$). However, the arguments in the last two arctangent function arguments are not integer reciprocals and, consequently, practical application of equation \eqref{eq_14} may be questionable. In order to eliminate these two quotients in arctangent function arguments we can use the following identity
$$
\arctan\left(\frac{1}{x}\right) = 2\arctan\left(\frac{1}{2x}\right) - \arctan\left(\frac{1}{4x^3 + 3x}\right)
$$
leading to
$$
\begin{aligned}
\arctan\left(\frac{2}{2513489}\right) = &2\arctan\left(\frac{1}{2513489}\right)\\
&-\arctan\left(\frac{1}{7939642926390344818}\right)
\end{aligned}
$$
and
$$
\begin{aligned}
\arctan\left(\frac{2}{18280007883}\right) = &2\arctan\left(\frac{1}{18280007883}\right)\\
&-\arctan\left(\frac{1}{3054211727257704725384731479018}\right).
\end{aligned}
$$
Substituting these two arctangent function relations into equation \eqref{eq_14} results in equation \eqref{eq_11}.

Consider now how our approach based on identity \eqref{eq_13} can be used to eliminate quotients in \eqref{eq_14} \cite{Abrarov2021a}. Specifically, applying equation \eqref{eq_13} in a sequence first for
$$
z=-\frac{2513489}{2}
$$
and then for
$$
z=-\frac{18280007883}{2},
$$
we obtain the following 7-term Machin-like formula for $\pi$
\begin{equation}\label{eq_15}
\begin{aligned}
\frac{\pi}{4} = &83\arctan\left(\frac{1}{107}\right) + 17\arctan\left(\frac{1}{1710}\right) - 22\arctan\left(\frac{1}{103697}\right) \\
&-12\arctan\left(\frac{1}{1256744}\right) - 22\arctan\left(\frac{1}{9140003941}\right) \\
&+12\arctan\left(\frac{1}{3158812219818}\right) \\
&+22\arctan\left(\frac{1}{167079344092131066905}\right),
\end{aligned}
\end{equation}
where all arguments are integer reciprocals. The Lehmer's measure for this equation is $\mu \approx 1.39524$. It seems that equation \eqref{eq_15} is new and first appeared in our recent publication \cite{Abrarov2021a}.

The Mathematica code below:\\[0.25cm]
\small
{\ttfamily coeff=(107+I)\string^83*(1710+I)\string^17*(103697+I)\string^-22*(1256744+I)\string^-12*

\noindent\quad(9140003941+I)\string^-22*(3158812219818+I)\string^12*

\noindent\qquad(167079344092131066905+I)\string^22;\\
Re[coeff]==Im[coeff]}\\[0.25cm]
\normalsize
validates equation \eqref{eq_15} by returning {\ttfamily True}.

In 2002, Kanada broke a record by computing more than one trillion decimal digits of $\pi$ \cite{Calcut2009,Agarwal2013} (current record is $62.8$ trillion digits of $\pi$ \cite{Thegardian2021}). He applied the following self-checking pair of the Machin-like formulas
\[
\begin{aligned}
\frac{\pi}{4} = &44\arctan\left(\frac{1}{57}\right) + 7\arctan\left(\frac{1}{239}\right) - 12\arctan\left(\frac{1}{682}\right)\\
&+24\arctan\left(\frac{1}{12943}\right)
\end{aligned}
\]
and
\[
\begin{aligned}
\frac{\pi}{4} = &12\arctan\left(\frac{1}{49}\right) + 32\arctan\left(\frac{1}{57}\right) - 5\arctan\left(\frac{1}{239}\right)\\
&+12\arctan\left(\frac{1}{110443}\right)
\end{aligned}
\]
with Lehmer's measures $1.58604$ and $1.7799$, respectively.

The current record that was achieved by using the Chudnovsky brothers formula \cite{Agarwal2013} is $62.8$ trillion digits of $\pi$ \cite{Thegardian2021}. This record is larger than Kanada's $20$ years old record by a factor $\sim 60$. Since this factor is not large, despite relatively high Lehmer's measures $1.58604$ and $1.7799$ we may expect that using most modern supercomputer the self-checking pair above can also generate a comparable amount of digits of $\pi$. Therefore, the Machin-like formulas with small Lehmer's measure have a colossal potential and remain competitive in computing digits of $\pi$.

Since Lehmer's measures of equations \eqref{eq_11} and \eqref{eq_15} are smaller, their application in computing $\pi$ is potentially more efficient than the self-checking pair that was used by Kanada. Furthermore, since first three arctangent function terms in both equations are same, it is easier to compare results of computation when equations \eqref{eq_11} and \eqref{eq_15} are used as a self-checking pair.

We have shown already that identity \eqref{eq_13} can be used in an iterative procedure such that from equation \eqref{eq_12} it follows that (see \cite{Abrarov2021a} for derivation)
\begin{align}\label{eq_16}
\frac{\pi}{4} = &8\arctan\frac{1}{10} - \arctan\frac{1}{84} - \arctan\frac{1}{21342} \nonumber \\
&- \arctan\frac{1}{991268848} - \arctan\frac{1}{193018008592515208050} \nonumber \\
&- \arctan\frac{1}{197967899896401851763240424238758988350338}\\
&- \arctan\scalebox{0.8}[1.3]{$\frac{1}{117573868168175352930277752844194126767991915008537018836932014293678271636885792397}$}, \nonumber
\end{align}
As we can see, this formula contains only integer reciprocals. Therefore, the iteration procedure based on equation \eqref{eq_13} can be used as an alternative to the Todd's process.

The 7-term Machin-like formula \eqref{eq_16} for $\pi$ can be validated by running the following Mathematica code: \\[0.25cm]
\small
{\ttfamily
coeff=(10+I)\string^8*(84+I)\string^-1*(21342+I)\string^-1*(991268848+I)\string^-1*

\noindent\quad(193018008592515208050+I)\string^-1*

\noindent\qquad(197967899896401851763240424238758988350338+I)\string^-1*

\noindent\qquad\quad(1175738681681753529302777528441941267679919150085370188369\textbackslash

\noindent\qquad\qquad 32014293678271636885792397+I)\string^-1;\\
Re[coeff]==Im[coeff]} \\[0.25cm]
\normalsize
The output of this code is {\ttfamily True}.

Thus, the initial equation \eqref{eq_12} with quotient can be transformed into equation \eqref{eq_16} with help of the identity \eqref{eq_13} using an iterative procedure described in our publication \cite{Abrarov2021a}. We can also notice that starting from the second arctangent function term each consecutive integer in equation \eqref{eq_16} is larger by absolute value than the previous one by many orders of the magnitude. This signifies that each consecutive arctangent term contributes smaller and smaller amount for the Lehmer measure $\mu$. Consequently, due to rapidly increasing integers this iteration process enable us to generate the Machin-like formula for $\pi$ with a small Lehmer measure.

\section{Derivation}

We will attempt to transform the Machin-like formula \eqref{eq_12} for $\pi$ in a form where initial arctangent terms contain arguments of kind $1/{10^m}$, where $m$ is a positive integer. Solving the equation
$$
\frac{1758719}{147153121} + x = \frac{1}{100}\Leftrightarrow x = -\frac{28718779}{14715312100}
$$
we can rewrite equation \eqref{eq_12} as follows
$$
\frac{\pi}{4} = 8\arctan\left(\frac{1}{10}\right) - \arctan\left(\frac{1}{100} + \frac{28718779}{14715312100}\right)
$$
Applying now the following identity (see also \cite{Abrarov2021a})
$$
\arctan\left(x + y\right) = \arctan\left(x\right) + \arctan\left(\frac{y}{1 + \left(x + y\right)x} \right)
$$
the equation \eqref{eq_12} can be recast as
\begin{equation}\label{eq_17}
\frac{\pi}{4} = 8\arctan\left(\frac{1}{10}\right) - \arctan\left(\frac{1}{100}\right) - \arctan\left( \frac{28718779}{14717070819}\right).
\end{equation}

Repeating same procedure again for equation \eqref{eq_17}
$$
\frac{28718779}{14717070819} + x = \frac{1}{1000}\Leftrightarrow x = -\frac{14001708181}{14717070819000}
$$
we get
\begin{equation}\label{eq_18}
\begin{aligned}
\frac{\pi}{4} = &8\arctan\left(\frac{1}{10}\right) - \arctan\left(\frac{1}{100}\right) \\
&-\arctan\left(\frac{1}{1000} + \frac{14001708181}{14717070819000}\right) \\ 
= &8\arctan \left(\frac{1}{10}\right) - \arctan\left(\frac{1}{100}\right) - \arctan\left(\frac{1}{1000}\right) \\
&-\arctan\left(\frac{14001708181}{14717099537779}\right).
\end{aligned}
\end{equation}

Since the following quotient
$$
\frac{14717099537779}{14001708181}\approx 1051.09
$$
is close to $1000$, it is reasonable to represent the last arctangent function from equation \eqref{eq_18} as
$$
\arctan\left(\frac{14001708181}{14717099537779}\right) = \arctan\left(\frac{1}{1000} - \frac{715391356779}{14717099537779000}\right).
$$
This leads to
\begin{equation}\label{eq_19}
\begin{aligned}
\frac{\pi}{4} = &8\arctan\left(\frac{1}{10}\right) - \arctan\left(\frac{1}{100}\right) - 2\arctan\left(\frac{1}{1000}\right) \\
&+\arctan\left(\frac{715391356779}{14717113539487181}\right).
\end{aligned}
\end{equation}

The following Mathematica code:\\[0.25cm]
\small
{\ttfamily
coeff=(10+I)\string^8*(100+I)\string^-1*

\noindent\quad(14717070819/28718779+I)\string^-1;\\
eq17=Re[coeff]==Im[coeff];\\[0.25cm]
coeff=(10+I)\string^8*(100+I)\string^-1*

\noindent\quad(14717070819/28718779+I)\string^-1;\\
eq18=Re[coeff]==Im[coeff];\\[0.25cm]
coeff=(10+I)\string^8*(100+I)\string^-1*(1000+I)\string^-1*

\noindent\quad(14717099537779/14001708181+I)\string^-1;\\
eq19=Re[coeff]==Im[coeff];\\[0.25cm]
Print[\{eq17,eq18,eq19\}]}\\[0.25cm]
\normalsize
validates equations \eqref{eq_17}, \eqref{eq_18} and \eqref{eq_19} by returning the list {\ttfamily\{True,True,True\}}.

We notice that
$$
\left\lfloor\frac{14717113539487181}{715391356779}\right\rfloor = 20572
$$
is a relatively large number. Therefore, its contribution to the Lehmer measure will not be high and  further we can apply identity \eqref{eq_13}. Thus, taking
$$
z = -\frac{14717113539487181}{715391356779}
$$
and substituting it into in identity \eqref{eq_13} we can rewrite equation \eqref{eq_19} as
\begin{equation}\label{eq_20}
\begin{aligned}
\frac{\pi}{4} = &8\arctan\left(\frac{1}{10}\right) - \arctan\left(\frac{1}{100}\right) - 2\arctan\left(\frac{1}{1000}\right) \\
&+\arctan\left(\frac{1}{20573}\right) + \arctan\left(\frac{316421763593}{151387588781630565746}\right).
\end{aligned}
\end{equation}

Now implying that
$$
z = -\frac{151387588781630565746}{316421763593}
$$
and substituting it into identity \eqref{eq_13}, we obtain
\begin{equation}\label{eq_21}
\begin{aligned}
\frac{\pi}{4} = &8\arctan\left(\frac{1}{10}\right) - \arctan\left(\frac{1}{100}\right) - 2\arctan\left(\frac{1}{1000}\right) \\
&+\arctan\left(\frac{1}{20573}\right) +\arctan\left(\frac{1}{478436082}\right) \\
&+\arctan\left(\frac{10266919376}{14485856968022096352676282153}\right).
\end{aligned}
\end{equation}

The Mathematica code below validates equations \eqref{eq_20} and \eqref{eq_21}:\\[0.25cm]
\small
{\ttfamily
coeff=(10+I)\string^8*(100+I)\string^-1*(1000+I)\string^-2*(20573+I)*

\noindent\quad(151387588781630565746/316421763593+I);\\[0.25cm]
eq20=Re[coeff]==Im[coeff];\\[0.25cm]
coeff=(10+I)\string^8*(100+I)\string^-1*(1000+I)\string^-2*

\noindent\quad(20573+I)*(478436082+I)*

\noindent\qquad(14485856968022096352676282153/10266919376+I);\\[0.25cm]
eq21=Re[coeff]==Im[coeff];\\[0.25cm]
Print[\{eq20,eq21\}]}\\[0.25cm]
\normalsize
\noindent by returning the list {\ttfamily \{True,True\}}.

Repeating same procedure over and over again, we end up with the following $21$-term Machin-like formula for $\pi$
\begin{equation}\label{eq_22}
\begin{aligned}
\frac{\pi}{4} = &8\arctan\left(\frac{1}{10}\right) - \arctan\left(\frac{1}{100}\right) - 2\arctan\left(\frac{1}{1000}\right) \\
&+\arctan\left(\frac{1}{20573}\right) + \arctan\left(\frac{1}{478436082}\right) \\
&+\arctan\left(\frac{1}{1410925365001336732}\right)\\
&+\arctan\left(\frac{1}{2921851992939769423775706369842706095}\right)\\ 
&+\cdots + \arctan \left(\frac{1}{\underbrace{1170619828 \ldots 6246893153}_{\text{600,593 digits}}}\right),
\end{aligned}
\end{equation}
where integer in the last term consists of $600,593$ digits.

It is not difficult to show that equation \eqref{eq_22} can be represented in a more compact form as
\[
\begin{aligned}
\frac{\pi}{4} =&8\arctan\left(\frac{1}{10}\right) - \arctan\left(\frac{1}{100}\right) - 2\arctan\left(\frac{1}{1000}\right) \\
&-\sum\limits_{m = 1}^{18}\arctan\left(\frac{1}{\left\lfloor{\cal B}_m\right\rfloor}\right),
\end{aligned}
\]
where
$$
{\cal B}_1 = -\frac{14717113539487181}{715391356779}
$$
such that by iteration
$$
{\cal B}_m = \frac{1 + \left\lfloor{\cal B}_{m - 1}\right\rfloor{\cal B}_{m - 1}}{\left\lfloor{\cal B}_{m - 1}\right\rfloor - {\cal B}_{m - 1}}.
$$
We ceased iteration at $18^{th}$ step since the number
$$
{\cal B}_{18} = \left\lfloor{\cal B}_{18}\right\rfloor = -\underbrace{1170619828\ldots 6246893153}_{\text{600,593 digits}}
$$
is not a quotient but an integer (see \cite{Abrarov2021a} showing how iteration is implemented).

The following Mathematica code:\\[0.25cm]
\small
{\ttfamily
Clear[B]\\[0.25cm]
B[1]:=B[1]=-14717113539487181/715391356779;\\
B[m\char95]:=B[m]=(1+Floor[B[m-1]]*B[m-1])/

\noindent\quad(Floor[B[m-1]]-B[m-1]);\\[0.25cm]
coeff=(10+I)\string^8*(100+I)\string^-1*

\noindent\quad(1000+I)\string^-2*Product[(Floor[B[m]]+I)\string^-1,

\noindent\qquad \{m,1,18\}];\\[0.25cm]
Re[coeff]==Im[coeff]}\\[0.25cm]
\normalsize
returns {\ttfamily True} for validation of the equation \eqref{eq_22}.

As we have reported in \cite{Abrarov2021b}, the arctangent function can be expanded as a series
$$
\arctan\left(x\right)=\sum_{m = 1}^M\arctan\left(\frac{M x}{M^2 + \left(m - 1\right)\,m x^2}\right).
$$
Using this series expansion we can obtain the following relation
$$
\begin{aligned}
\arctan\left(\frac{1}{10}\right) = &10\arctan\left(\frac{1}{100}\right) \\
&\hspace{-1.25cm}-\frac{1}{2}\arctan\left(\frac{663251878199233259534207841462897801980}{1011010237065062789681261499707188111860099}\right),
\end{aligned}
$$
where the quotient on the right side can be eliminated by using equation \eqref{eq_13} in iterative procedure. Therefore, there is a possibility to rearrange the equation \eqref{eq_12} with starting term $80\arctan(1/100)$. This may be more advantageous in computation due to faster convergence.

\section{Lehmer's measure}

According to definition \eqref{eq_3}, the Lehmer measure for the equation \eqref{eq_22} is
$$
\frac{1}{\log_{10}\left(10\right)} + \frac{1}{\log_{10}\left(100\right)} + \frac{1}{\log_{10}\left(1000 \right)} + \sum\limits_{j = 1}^{18}\frac{1}{\log_{10}\left(\left|\left\lfloor{\cal B}_j\right\rfloor\right|\right)} \approx 2.29025.
$$
However, the measure may be reduced under some specific criteria. In particular, Lehmer justifiably noticed that the actangent terms of kind
$$
\arctan\left(\frac{1}{{10}^m}\right),\quad m\in\mathbb{N}.
$$
are much easier to calculate (see \cite{Lehmer1938}, page 662). Therefore, he suggested to assign the measure $1/2$ for the arctangent term $\arctan\left(1/10\right)$. Moreover, Lehmer also pointed out that if a Machin-like formula for $\pi$ besides $\arctan\left(1/10\right)$ includes the arctangents of other powers of 1/10, then the measure $0$ should be assigned to all of them since the Maclaurin series expansion \eqref{eq_2} of $\arctan \left(1/10\right)$ already contains all required digits in mantissa of each term that can be used further for computation of $\arctan\left(1/10^m\right)$. More explicitly, once in accordance with equation \eqref{eq_2} we compute
$$
\arctan\left(\frac{1}{10}\right) = \frac{1}{10} - \frac{1}{3\cdot 10^3} + \frac{1}{5\cdot 10^5} - \frac{1}{7\cdot 10^7}\dots\, ,
$$
then the computation
\small
$$
\begin{aligned}
\arctan\left(\frac{1}{10^m}\right) = &\left(\frac{1}{10}\right)\frac{1}{10^{m - 1}} + \left(-\frac{1}{3\cdot 10^3}\right)\frac{1}{10^{3\left(m - 1\right)}}\\
&+\left(\frac{1}{5\cdot 10^5}\right)\frac{1}{10^{5\left(m - 1\right)}} + \left(-\frac{1}{7\cdot 10^7}\right)\frac{1}{10^{7\left(m - 1\right)}}\dots
\end{aligned}
$$
\normalsize
becomes straightforward since it can be implemented algorithmically just by decreasing exponent in each term while keeping mantissa absolutely unchanged. Consequently, as Lehmer suggested \cite{Lehmer1938} if $\arctan\left(1/10\right)$ is present in a given relation, then for each arctangent term $\arctan\left(1/{10}^m\right)$ the measure $0$ can be assigned. Thus, the Lehmer measure for the equation \eqref{eq_22} can be reduced as
$$
\underbrace{\frac{1}{2}}_{\frac{1}{\log_{10}\left(10\right)}} + \underbrace{0}_{\frac{1}{\log_{10}\left(100\right)}} + \underbrace{0}_{\frac{1}{\log_{10}\left(1000\right)}} + \sum\limits_{j = 1}^{18}\frac{1}{\log_{10}\left(\left|\left\lfloor{\cal B}_j\right\rfloor\right|\right)} \approx 0.956915.
$$

The corresponding Mathematica code is:\\[0.25cm]
\small
{\ttfamily
Clear[B]\\[0.25cm]
B[1]:=B[1]=-14717113539487181/715391356779;\\
B[m\char95]:=B[m]=(1+Floor[B[m-1]]*B[m-1])/

\noindent\quad(Floor[B[m-1]]-B[m-1]);\\[0.25cm]
Print["The Lehmer measure of identity (22) is: ",

\noindent\quad 1/2+0+0+Sum[1/Log10[Abs[B[m]]],

\noindent\qquad\quad \{m,1,18\}]//N]\\[0.25cm]
}
\normalsize
The generated output is:\\[0.25cm]
{\ttfamily
The Lehmer measure of identity (22) is: 0.956916\\[0.25cm]
}

As we can see from this example, the Lehmer measure below one can be achieved in a Machin-like formula for $\pi$ without any quotient and at a smaller values of the integer $k$. Despite that the initial formula \eqref{eq_12} for $\pi$ that we used as a starting point in iterative procedure, corresponds to the small integer $k = 4$, this method enables us to generate the multi-term Machin-like formula \eqref{eq_22} for $\pi$ with Lehmer measure less than unity, $\mu \approx 0.956915$.

Consider for comparison the following two Machin-like formulas for $\pi$ that also include three arctangent terms with arguments $1/10$, $1/100$ and $1/1000$ \cite{Wrench1938}
\begin{equation}\label{eq_23}
\begin{aligned}
\frac{\pi}{4} = &7\arctan\left(\frac{1}{10}\right) + 2\arctan\left(\frac{1}{50}\right) + 4\arctan\left( \frac{1}{100}\right) \\
&+\arctan\left(\frac{1}{682}\right) + 4\arctan\left(\frac{1}{1000}\right) + 3\arctan\left(\frac{1}{1303}\right) \\
&-4\arctan\left(\frac{1}{90109}\right)
\end{aligned}
\end{equation}
and
\begin{equation}\label{eq_24}
\begin{aligned}
\frac{\pi}{4} = &7\arctan \left(\frac{1}{10}\right) + 8\arctan\left(\frac{1}{100}\right) + \arctan\left(\frac{1}{682}\right) \\
&+4\arctan\left(\frac{1}{1000}\right) + 3\arctan\left(\frac{1}{1303}\right) - 4\arctan\left(\frac{1}{90109}\right) \\
&-2\arctan\left(\frac{1}{500150}\right).
\end{aligned} 
\end{equation}
The corresponding Lehmer measures are
\[
\begin{aligned}
&\frac{1}{2} + \frac{1}{\log_{10}\left(50\right)} + 0 + \frac{1}{\log_{10}(682)} + 0 + \frac{1}{\log_{10}\left(1303\right)} \\
&+ \frac{1}{\log_{10}\left(90109\right)} \approx 1.96434
\end{aligned}
\]
and
\[
\begin{aligned}
&\frac{1}{2} + 0 + \frac{1}{\log_{10}\left(682\right)} + 0 + \frac{1}{\log_{10}\left(1303\right)} + \frac{1}{\log_{10}\left(90109\right)} \\
&+ \frac{1}{\log_{10}\left(500150\right)} \approx 1.55121,
\end{aligned}
\]
respectively. As we can see, the Machin-like formula \eqref{eq_22} for $\pi$ is computationally more efficient since its Lehmer's measure is significantly smaller.

The computation of the first three arctangent terms in equation \eqref{eq_22} can be implemented by using the Maclaurin series expansion \eqref{eq_2} since in decimal system no actual multiplications in determination of powers of 1/10 are needed (as mentioned, this can be done only by decreasing exponent and keeping mantissa unchanged in each term). However, starting from the forth arctangent term $\arctan\left(1/20573\right)$ application of the Maclaurin series expansion \eqref{eq_2} may not be optimal since faster convergence rate can be achieved by using either Euler's formula \cite{Castellanos1988,Chien-Lih2005}
\begin{equation}\label{eq_25}
\arctan\left(x\right) = \sum\limits_{m = 0}^\infty\frac{2^{2m}\left(m!\right)^2}{\left(2m + 1\right)!}\frac{x^{2m + 1}}{\left(1 + x^2\right)^{m + 1}}
\end{equation}
or iteration-based formula \cite{Abrarov2021b}
\begin{equation}\label{eq_26}
\arctan\left(x\right) = 2\sum\limits_{m = 1}^\infty  \frac{1}{2m - 1}\frac{g_m\left(x\right)}{g_m^2\left(x\right) + h_m^2\left(x\right)},
\end{equation}
where
$$
g_1\left(x\right) = 2/x,h_1\left(x\right) = 1,
$$
$$
g_m\left(x\right) = g_{m - 1}\left(x\right)\left(1 - 4/x^2\right) + 4h_{m - 1}\left(x\right)/x,
$$
$$
h_m\left(x\right) = h_{m - 1}\left(x\right)\left(1 - 4/x^2\right) - 4g_{m - 1}\left(x\right)/x.
$$

Chien-Lih showed an elegant derivation of the Euler's formula \eqref{eq_25} by an elementary method \cite{Chien-Lih2005}. In particular, this series expansion can be derived by taking integral
$$
\arctan\left(x\right) = \int_0^{\pi/2}\frac{x\sin\left(u\right)}{1+x^2}\;\frac{1}{\left(1 - \frac{x^2\sin^2\left(u\right)}{1 + x^2}\right)}du
$$
in terms of geometric series
$$
\frac{1}{1 - \frac{x^2\sin^2\left(u\right)}{1 + x^2}} = \sum_{n=0}^{\infty}\frac{x^{2n}\sin^{2n}\left(u\right)}{\left(1 + x^2\right)^n}.
$$

Iteration-based equation \eqref{eq_26} represents a trivial rearrangement of the series expansion
\[
\arctan\left(x\right) = i\sum_{n = 1}^{\infty}\frac{1}{2n - 1}\left(\frac{1}{\left(1 + 2i/x\right)^{2n - 1}} - \frac{1}{\left(1 - 2i/x\right)^{2n - 1}}\right)
\]
that we derived in our work \cite{Abrarov2018}.

It should be noted that apart from decimal system, other number systems can also be used in the Machin-like formulas for $\pi$. For example, a more rapid computation of the Maclaurin series expansion \eqref{eq_2} can also be achieved for the arctangent terms like $\arctan\left(1/{8^m}\right)$ and $\arctan\left(1/{16}^m\right)$ in octal and hexadecimal systems, respectively.

\section{Klingenstierna's identity}

The method described above can also be used to derive some known Machin-like formulas for $\pi$. This derivation is much easier than the conventional derivations. For example, solving
$$
\frac{28718779}{14717070819} + x = \frac{1}{515}\Leftrightarrow x = -\frac{73100366}{7579291471785}
$$
and representing the last arctangent term from equation \eqref{eq_17} as
$$
-\arctan\left(\frac{28718779}{14717070819}\right) = \arctan\left(-\frac{1}{515} - \frac{73100366}{7579291471785}\right)
$$
we immediately get the well-known identity
\[
\frac{\pi}{4} = 8\arctan\left(\frac{1}{10}\right) - \arctan\left(\frac{1}{100}\right) - \arctan\left(\frac{1}{515}\right) - \arctan\left(\frac{3583}{371498882}\right)
\]
According to Castellanos, this identity is attributed to Swedish mathematician Samuel Klingenstierna \cite{Castellanos1988}. Although the corresponding measure
\[
\frac{1}{2} + 0 + \frac{1}{\log_{10}\left(515\right)} + \frac{1}{\log_{10}\left(\frac{371498882}{3583}\right)} \approx 1.06813
\]
is relatively small, this Machin-like formula may not be optimal for computing $\pi$ due to quotient in its last arctangent term.

\section{Alternative representation}

It is interesting to note that a slight modification of the identity \eqref{eq_13} as given by
\begin{equation}\label{eq_27}
\begin{aligned}
\arctan\left(\frac{1}{z}\right) = &\arctan\left(\frac{1}{\left\lfloor {10}^n z\right\rfloor {10}^{-n}} \right) \\
&+\arctan\left(\frac{\left\lfloor{10}^n z \right\rfloor {10}^{-n} - z}{1 + z\left\lfloor{10}^n z \right\rfloor{10}^{-n}}\right),\qquad z \notin [0,10^{-n}),
\end{aligned}
\end{equation}
can also be used to generate the Machin-like formulas for $\pi$. Return to equation \eqref{eq_12} again and assume that
$$
z = -\frac{147153121}{1758719}
$$
For, say, $n = 2$ substituting this value $z$ into identity \eqref{eq_27} leads to
\begin{equation}\label{eq_28}
\frac{\pi}{4} = 8\arctan\left(\frac{1}{10}\right) - \arctan\left(\frac{1}{83.68}\right) - \arctan\left( \frac{412123}{307888297107}\right).
\end{equation}
Similarly, substituting now
$$
z = -\frac{307888297107}{412123}
$$
into identity \eqref{eq_27} results in
\begin{equation}\label{eq_29}
\begin{aligned}
\frac{\pi}{4} = &8\arctan \left(\frac{1}{10}\right) - \arctan\left(\frac{1}{83.68}\right) - \arctan\left(\frac{1}{747078.66}\right) \\
&-\arctan\left(\frac{74409}{11500838821639577981}\right).
\end{aligned}
\end{equation}
The following Mathematica code:\\[0.25cm]
\small
{\ttfamily
coeff=(10+I)\string^8*(83+68/100+I)\string^-1*

\noindent\quad(307888297107/412123+I)\string^-1;\\
eq28=Re[coeff]==Im[coeff];\\[0.25cm]
coeff=(10+I)\string^8*(83+68/100+I)\string^-1*(747078+66/100+I)\string^-1*

\noindent\quad(11500838821639577981/74409+I)\string^-1;\\
eq29=Re[coeff]==Im[coeff];\\[0.25cm]
Print[\{eq28,eq29\}]}\\[0.25cm]
validates equations \eqref{eq_28} and \eqref{eq_29} by returning the list {\ttfamily \{True,True\}}.

Using the iteration, we can derive the following Machin-like formula for $\pi$
\begin{equation}\label{eq_30}
\begin{aligned}
\frac{\pi}{4} = &8\arctan\left(\frac{1}{10}\right) - \arctan\left(\frac{1}{83.68}\right) - \arctan\left(\frac{1}{747078.66}\right) \\
&-\arctan\left(\frac{1}{154562469884551.31}\right) \\
&-\arctan \left(\frac{1}{399648835184411935214717088966.73}\right) \\
&-\cdots -\arctan\left(\frac{1}{\underbrace{8665971818\ldots 7222871549}_{\text{1,052 digits}}}\right) 
\end{aligned}
\end{equation}
consisting of 10-terms with finite number of decimal digits in each arctangent term.

The 10-term Machin-like formula \eqref{eq_30} for $\pi$ can be represented in a compact form as follows
$$
\frac{\pi}{4} = 8\arctan\left(\frac{1}{10}\right) + \sum\limits_{m = 1}^9\arctan\left(\frac{1}{\left\lfloor{10}^2{\cal B}_m\right\rfloor{10}^{-2}}\right),
$$
where
$$
{\cal B}_1 = -\frac{147153121}{1758719}
$$
and
$$
{\cal B}_m = \frac{1 + \left(\lfloor 10^2 {\cal B}_{m-1}\rfloor 10^{-2}\right){\cal B}_{m - 1}}{\left(\lfloor 10^2{\cal B}_{m - 1}\rfloor 10^{-2}\right)-{\cal B}_{m - 1}}.
$$
Iteration is completed at $9^{th}$ step since (see \cite{Abrarov2021a} for a detailed description of iteration procedure)
$$
{\cal B}_9=\lfloor{10^2\cal B}_9\rfloor 10^{-2} = -{\underbrace{8665971818 \ldots 7222871549}_\text{1,052 digits}}
$$
is an integer.

The following Mathematica program:\\[0.25cm]
\small
{\ttfamily
Clear[B]\\[0.25cm]
B[1]:=B[1]=-147153121/1758719;\\
B[m\char95]:=B[m]=(1+(Floor[10\string^2*B[m-1]]*10\string^-2)*

\noindent\quad B[m-1])/((Floor[10\string^2*B[m-1]]*

\noindent\qquad 10\string^-2)-B[m - 1]);\\[0.25cm]
coeff=(10+I)\string^8*

\noindent\quad Product[(Abs[Floor[10\string^2*B[m]]*

\noindent\qquad 10\string^-2]+I)\string^-1,\{m,1,9\}];\\[0.25cm]
Re[coeff]==Im[coeff]}\\[0.25cm]
validates equation \eqref{eq_30} by returning {\ttfamily \{True\}}.

\section{Conclusion}

A method that can be used to generate a multi-term Machin-like formula for $\pi$ with a reduced Lehmer measure is presented. We show that the Lehmer measure below $1$ can be achieved even at small value of the integer $k=4$. Specifically, as an example we derive the Machin-like formula \eqref{eq_22} for $\pi$ with the Lehmer measure $\mu \approx 0.956916$. This value is significantly smaller than the smallest Lehmer measure $1.34085$ known so far among the Machin-like formulas for $\pi$ consisting of only integer reciprocals.

\section*{Acknowledgment}
\addcontentsline{toc}{section}{Acknowledgment}
This work is supported by National Research Council Canada, Thoth Technology Inc., York University and Epic College of Technology.

\bigskip


\begin{thebibliography}{99}

\bibitem[1]{Lehmer1938}
\bibinfo{author}{Lehmer, D.H.} \bibinfo{year}{(1938)}
\newblock \bibinfo{title}{On Arccotangent Relations for $\pi$},
\newblock \bibinfo{journal}{American Mathematical Monthly}, \bibinfo{volume}{{\bf 45}}(10), \bibinfo{pages}{657--664}. \bibinfo{doi}{\url{https://doi.org/10.2307/2302434}}

\bibitem[2]{Abeles1993}
\bibinfo{author}{Abeles, F.F.} \bibinfo{year}{(1993)} \bibinfo{title}{Charles L. Dodgson's Geometric Approach to Arctangent Relations for Pi}, {{\bf 20}}(2), \bibinfo{pages}{151--159}. \bibinfo{doi}{\url{https://doi.org/10.1006/hmat.1993.1013}}

\bibitem[3]{Beckmann1971}
\bibinfo{author}{Beckmann, P.} \bibinfo{year}{(1971)} \bibinfo{title}{A History of pi}, \bibinfo{publisher}{Golem Press}, New York.

\bibitem[4]{Berggren2004}
\bibinfo{author}{Berggren, L., Borwein J. and Borwein P.} \bibinfo{year}{(2004)} \bibinfo{title}{Pi: a Source Book}, $\rm{3}^{\rm{rd}}$ ed., \bibinfo{publisher}{Springer-Verlag}, New York.

\bibitem[5]{Borwein2008}
\bibinfo{author}{Borwein J. and Bailey, D.} \bibinfo{year}{(2008)} \bibinfo{title}{Mathematics by Experiment. Plausible Reasoning in the $\rm{21^{st}}$ Century}, $\rm{2}^{\rm{nd}}$ ed., \bibinfo{publisher}{Taylor \& Francis Group}.

\bibitem[6]{Fitzhugh2013}
\bibinfo{author}{Fitzhugh, J.M. and Farnsworth, D.L.} \bibinfo{year}(2013) \bibinfo{title}{A Construction That Produces Wallis-Type Formulas}, \bibinfo{journal}{Advances in Pure Mathematics}, \bibinfo{volume}{{\bf 3}}, \bibinfo{pages}{579--585}. \bibinfo{doi}{\url{https://dx.doi.org/10.4236/apm.2013.36074}}

\bibitem[7]{Hwang2017}
\bibinfo{author}{Hwang, C.-O., Kim, Y., Im, C. and Lee S.} \bibinfo{year}{(2017)} \bibinfo{title}{Buffon's Needle Algorithm to Estimate $\pi$}, \bibinfo{journal}{Applied Mathematics}, \bibinfo{volume}{{\bf 8}}, \bibinfo{pages}{275--279}. \bibinfo{doi}{\url{https://doi.org/10.4236/am.2017.83022}}

\bibitem[8]{Brent2010}
\bibinfo{author}{Brent R.P. and Zimmermann. P.} \bibinfo{year}{(2010)} \bibinfo{title}{Modern Computer Arithmetic}, \bibinfo{publisher}{Cambridge University Press}, Cambridge.

\bibitem[9]{Tweddle1991}
\bibinfo{author}{Tweddle, I.} \bibinfo{year}{(1991)} \bibinfo{title}{John Machin and Robert Simson on Inverse-Tangent Series for $\pi$}, \bibinfo{journal}{Archive for History of Exact Sciences}, \bibinfo{volume}{{\bf 42}(1)}, \bibinfo{pages}{1--14}. \bibinfo{url}{\url{https://www.jstor.org/stable/41133896}}

\bibitem[10]{Wetherfield1996}
\bibinfo{author}{Wetherfield, M.} \bibinfo{year}{(1996)} \bibinfo{title}{The Enhancement of Machin's Formula by Todd's Process}, \bibinfo{journal}{Mathematical Gazette}, \bibinfo{volume}{{\bf 80}(488)}, \bibinfo{pages}{333--344}. \bibinfo{url}{\url{https://www.jstor.org/stable/3619567}}

\bibitem[11]{Jansson2019}
\bibinfo{author} Jansson, M. \bibinfo{year}{(2019)} \bibinfo{title}{\href{https://lup.lub.lu.se/student-papers/search/publication/8983341}{Approximation of $\pi$}}, \bibinfo{publisher}{Lund University}.

\bibitem[12]{Nimbran2010}
\bibinfo{author}{Nimbran, A.S.} \bibinfo{year}{(1990)} \bibinfo{title}{On the Derivation of Machin-like Arctangent Identities for Computing pi ($\pi$)}, \bibinfo{journal}{The Mathematics Student}, \bibinfo{volume}{{\bf 79}(1--4)}, \bibinfo{pages}{171--186}.

\bibitem[13]{Chamberland2019}
\bibinfo{author}{Chamberland M. and Herman, E.A.} \bibinfo{year}{(2019)} \bibinfo{title}{Arctangent Formulas and pi}, \bibinfo{journal}{American Mathematical Monthly}, \bibinfo{volume}{{\bf 126}(7)} \bibinfo{pages}{646--650}. \bibinfo{doi}{\url{https://doi.org/10.1080/00029890.2019.1606578}}

\bibitem[14]{Abrarov2017a}
\bibinfo{author}{Abrarov, S.M. and B.M. Quine} \bibinfo{year}{(2017)} \bibinfo{title}{The Two-term Machin-like Formula for pi with Small Arguments of the Arctangent Function}, \bibinfo{preprint}{arXiv:1704.02875}. \bibinfo{doi}{\url{https://doi.org/10.48550/arXiv.1704.02875}}

\bibitem[15]{WolframCloud}
\bibinfo{title}{A Wolfram Notebook Playing with Machin-Like Formulas}. \bibinfo{url}{\url{https://www.wolframcloud.com/obj/exploration/MachinLike.nb}}

\bibitem[16]{Abrarov2017b}
\bibinfo{author}{Abrarov, S.M. and Quine, B.M.} \bibinfo{year}{(2017)} \bibinfo{title}{An Iteration Procedure for a Two-term Machin-like Formula for pi with Small Lehmer's Measure}, \bibinfo{preprint}{arXiv:1706.08835}. \bibinfo{doi}{\url{https://doi.org/10.48550/arXiv.1706.08835}}

\bibitem[17]{Todd1949}
\bibinfo{author}{Todd, J.} \bibinfo{year}{(1949)} \bibinfo{title}{Problem on Arc Tangent Relations}, \bibinfo{American Mathematical Monthly}, \bibinfo{volume}{{\bf 56}(8)}, \bibinfo{pages}{517--528}. \bibinfo{doi}{\url{https://doi.org/10.1080/00029890.1949.11999434}}

\bibitem[18]{Wetherfield2004}
\bibinfo{title}{Identity Lists - Classification, Layout and Links}. \bibinfo{url}{\url{http://www.machination.eclipse.co.uk/IdLists.html}}

\bibitem[19]{Abrarov2021a}
\bibinfo{author}{Abrarov, S.M., Siddiqui, R., Jagpal, R.K. and Quine, B.M.}, \bibinfo{year}{(2022)} \bibinfo{title}{A New Form of the Machin-like Formula for pi by Iteration with Increasing Integers}, \bibinfo{journal}{Journal of Integer Sequences}, \bibinfo{volume}{{\bf 25}}, \bibinfo{item}{22.4.5}. \bibinfo{url}{\url{https://cs.uwaterloo.ca/journals/JIS/VOL25/Abrarov/abrarov5.pdf}}

\bibitem[20]{Guillera2009}
\bibinfo{author}{Guillera, J.} \bibinfo{year}{(2009)} \bibinfo{title}{History of the Formulas and Algorithms for pi}, \bibinfo{preprint}{arXiv:0807.0872}. \bibinfo{doi}{\url{https://doi.org/10.48550/arXiv.0807.0872}}

\bibitem[21]{Calcut2009}
\bibinfo{author}{Calcut, J.S.} \bibinfo{year}{(2009)} \bibinfo{title}{Gaussian Integers and Arctangent Identities for $\pi$}, \bibinfo{journal}{American Mathemtical Monthly}, \bibinfo{volume}{{\bf 116}(6)}\bibinfo{pages}{515--530}. \bibinfo{url}{\url{www.jstor.org/stable/40391144}}

\bibitem[22]{Agarwal2013}
\bibinfo{author}{Agarwal, R.P., Agarwal, H. and Sen, S.K.} \bibinfo{year}{(2013)} \bibinfo{title}{Birth, Growth and Computation of pi to Ten Trillion Digits}, \bibinfo{journal}{Advances in Difference Equations}, \bibinfo{item}{100}. \bibinfo{doi}{\url{https://doi.org/10.1186/1687-1847-2013-100}}

\bibitem[23]{Thegardian2021}
\bibinfo{title}{Swiss Researchers Calculate pi to New Record of 62.8tn Figures}, \bibinfo{publisher}{The Guardian}. \bibinfo{url}{\url{https://tinyurl.com/6e96nudh}}

\bibitem[24]{Abrarov2021b}
\bibinfo{author}{Abrarov, S.M., Siddiqui, R., Jagpal, R.K. and Quine, B.M.} \bibinfo{year}{(2021)} \bibinfo{title}{Unconditional Applicability of Lehmer's Measure to the Two-term Machin-like Formula for $\pi$}, \bibinfo{journal}{The Mathematica Journal}, \bibinfo{volume}{{\bf 23}}, \bibinfo{item}{2}. \bibinfo{doi}{\url{https://doi.org/10.3888/tmj.23-2}}

\bibitem[25]{Wrench1938}
\bibinfo{author}{Wrench, J.W. Jr.}, \bibinfo{year}{(1938)} \bibinfo{title}{On the Derivation of Arctangent Equalities}, \bibinfo{journal}{American Mathematical Monthly}, \bibinfo{volume}{{\bf 45}(2)}, \bibinfo{pages}{108--109}. \bibinfo{doi}{\url{https://doi.org/10.2307/2304280}}

\bibitem[26]{Castellanos1988}
\bibinfo{author}{Castellanos, D.} \bibinfo{year}{(1988)} \bibinfo{title}{The Ubiquitous $\pi$}, \bibinfo{journal}{Mathematics Magazine}, \bibinfo{volume}{{\bf 61}(2)}, \bibinfo{pages}{67--98}. \bibinfo{doi}{\url{https://doi.org/10.2307/2690037}}

\bibitem[27]{Chien-Lih2005}
\bibinfo{author}{Chien-Lih, H.} \bibinfo{year}{(2005)} \bibinfo{title}{An Elementary Derivation of Euler's Series for the Arctangent Function}, \bibinfo{volume}{{\bf 89}(516)}, \bibinfo{pages}{469--470}. \bibinfo{doi}{\url{https://doi.org/10.1017/S0025557200178404}}

\bibitem[28]{Abrarov2018}
\bibinfo{author}{Abrarov S.M. and Quine B.M.} \bibinfo{year}{(2018)} \bibinfo{title}{A Formula for pi Involving Nested Radicals}, \bibinfo{journal}{The Ramanujan Journal}, \bibinfo{volume}{{\bf 46}(3)}, \bibinfo{pages}{657--665}. \bibinfo{doi}{\url{https://doi.org/10.1007/s11139-018-9996-8}}

\end{thebibliography}
\end{document}